\documentclass{amsart}

\usepackage[mathscr]{eucal}
\usepackage{amsfonts}
\usepackage{amsmath}
\usepackage{amsthm}
\usepackage{amssymb}
\usepackage{latexsym}
\usepackage{epsfig}
\usepackage{hhline}

\newtheorem{theorem}{Theorem}[subsection]

\newtheorem{proposition}[theorem]{Proposition}

\theoremstyle{definition}
\newtheorem{definition}[theorem]{Definition}

\newtheorem{remark}[theorem]{Remark}

\newcommand{\Tr}{\text{Tr}}

\newcommand{\Hom}{\text{Hom}}

\newcommand{\linob}{\text{span}}

\newcommand{\eps}{\varepsilon}

\renewcommand{\O}{\mathcal{O}}

\newcommand{\h}{\mathfrak{h}}

\newcommand{\half}{1/2}

\newcommand{\ot}{\otimes}

\newcommand{\ben}{\begin{enumerate}}
\newcommand{\een}{\end{enumerate}}
\newcommand{\beq}{\begin{equation}}
\newcommand{\eeq}{\end{equation}}

\newcommand{\CC}{{\mathbb{C}}}

\newcommand{\p}{\partial}

\begin{document}

\title[Representations of the rational
Cherednik algebras]{Representations of the rational Cherednik
algebras of  dihedral type}

\begin{abstract} In this paper  we describe the Jordan-H{\"o}lder series of
the standard modules over the rational Cherednik algebras
associated with the dihedral group. In particular, we compute the
characters of the irreducible representations from the category
$\mathcal O$ and classify the finite dimensional representations.
\end{abstract}

\author{Tatyana Chmutova}
\address{Department of Mathematics, Harvard University,
Cambridge, MA 02138, USA}
\email{chmutova@math.harvard.edu}

\maketitle
\section{Introduction}
The rational Cherednik algebra $H_{\sf k}(W)$ for a complex
reflection group $W$ arises from the study of Dunkl operators
(\cite{Du},\cite{DO}), and is an important special case of
symplectic reflection algebras (\cite{EG}).

In the recent years, a number of authors studied representations
of $H_{\sf k}(W)$ (see e.g. \cite{GGOR,BEG,BEG1,DO, De, Go}).
Specifically, the paper \cite{GGOR} defines the category $\mathcal
O$ for $H_{\sf k}(W)$, and shows that its basic structure is quite
similar to the structure of the category $\mathcal O$ for a simple
Lie algebra, studied by Bernstein, I. Gelfand and S. Gelfand. In
particular, the main problem about this category is to find the
characters of the irreducible modules in $\mathcal O$ and their
multiplicities in standard modules.

Although many partial results are available, this problem remains
open in all interesting situations. The goal of this paper is to
solve this problem for dihedral groups, which is the simplest
nontrivial case. We describe the structure of the standard modules
over rational Cherednik algebras corresponding to the dihedral
group $I_2(d)$. The structure depends on the parity of $d$.

When $d$ is odd the Cherednik algebra depends on one complex
parameter. For general value of the parameter, all the standard
modules are irreducible. The set of exceptional values consists of
a finite number of arithmetic progressions of rational numbers,
each with the common difference $1$. For these values, the
standard modules have length at most $2$, and at least one of them
is reducible. The category $\mathcal O$ is not semisimple in this
case. This classification is given in theorems \ref{odd1dim} and
\ref{odd2dim}.

The case of even $d$ is more complicated since the Cherednik
algebra depends on two complex parameters. If the values of the
parameters do not belong to four exceptional families of lines,
all the standard modules are irreducible. For generic points of
the exceptional lines the maximal length of standard modules is
$2$ and at least one is reducible. The category $\mathcal O$ for
these values of parameters is not semisimple. At the points of
intersection of the lines in these families, the structure of the
standard modules might be complicated, and the length might reach
$4$. For some of the intersection points even the category of the
finite dimensional modules is not semisimple. This case is
described in theorems \ref{even1dim} and \ref{even2dim}.

\section{Preliminaries}
\subsection{}
Let $W$ be a finite real reflection group, $\h_{\mathbb R}$ be its
reflection representation, and $\h$ be the complexification of
this representation. Fix a non-degenerate $W$-invariant scalar
product $(-,-)$ on $\h$. Denote by the same sign the induced
scalar product on $\h^*$. Let $S \subset W$ be the set of all
reflections $s \in W$. To each reflection $s \in S$ we associate
its fixed hyperplane $H_s$, a nonzero linear function $\alpha_s
\in \h^*$ which vanishes on $H_s$, and the element
$\alpha_s^{\vee}=\frac{2(\alpha_s, - )}{(\alpha_s, \alpha_s)} \in
\h$. In other words, to each $s$ we associate the corresponding
positive root $\alpha_s$  and coroot $\alpha_s^{\vee}$.

Let us fix a $W$-invariant function ${\sf k}: S \to \CC$, where
$W$ acts on $S$ via conjugation.

\begin{definition}
{\it The rational Cherednik algebra} $H_{\sf k}(W)$ is an
associative algebra generated by all elements of $W$, $\h^*$, and
$\h$ with the following defining relations:

$$
wxw^{-1}=w(x), \quad wyw^{-1}=w(y), \quad [x, x']=0, \quad [y,
y']=0,
$$
$$
[y,x]=(y,x)-\sum_{s \in S}{\sf k}(s)(y, \alpha_s)(\alpha_s^{\vee},
x)s,
$$
where $x, x' \in \h^*$, $y, y' \in \h$ and $w \in W$ (see also
\cite{EG}).
\end{definition}

\subsection{} \label{standmod}
For rational Cherednik algebras there exists an analog of Verma
modules (see \cite{BEG}), called {\it standard} modules. They are
indexed by irreducible representations of $W$. Let $V$ be an
irreducible representation of $W$. It can be considered as a
representation of $\CC[W] \ltimes \CC[\h^*]$, where $y \in \h$
acts by $0$.

\begin{definition} {\it The standard module} $M_{\sf k}(V)$ over $H_{\sf k}(W)$
corresponding to $V$ is defined by the formula
$$
M_{\sf k}(V)= H_{\sf k}(W)\ot_{\CC[W] \ltimes \CC[\h^*]} V.
$$
\end{definition}

\begin{definition} The standard module
$M_{\sf k}(V)$ and its quotients are called {\it modules with
lowest $W$-weight $V$.}
\end{definition}

As a vector space,  $M_{\sf k}(V)$ is naturally identified with
$\CC[\h] \ot V$. Hence the usual grading on $\CC[\h]$ induces a
grading on $M_{\sf k} (V)$.

The standard module $M_{\sf k}(V)$ has a unique simple quotient
$L_{\sf k}(V)$ (see \cite{DO}).

An important special case of standard modules is $M_{\sf k}({\sf
triv})$, which is also called the {\it polynomial representation}
of $H_{\sf k}(W)$. The space of this representation is $\CC [\h]$.
An element $w \in W$ acts on $\CC [\h]$ in the natural way, $x\in
\h$ acts by left multiplication, and the action of $y \in \h$ is
given by the {\it Dunkl operator}

\beq \label{Dunkl}
 D_y=\p_y-\sum_{s \in S}{\sf k}(s)\frac{(y,\alpha_s)}{\alpha_s} (1-s).
 \eeq

\begin{definition}
A vector in $M_{\sf k}(V)$ is called {\it singular} if it is
annihilated by every $y \in \h$.
\end{definition}

In particular, the singular vectors in $M_{\sf k}({\sf triv})$ are
exactly those annihilated by all Dunkl operators.

\subsection{} \label{O}
Every rational Cherednik algebra contains a canonically defined
$\mathfrak{sl}_2$-triple:

$$
{\mathbf h}=\frac{1}{2}\sum_i(x_iy_i+ y_ix_i), \quad E=\frac{1}{2}
\sum_i x_i^2, \quad F=-\frac{1}{2} \sum_i y_i^2,
$$
where $\{x_i\}$ and $\{y_i\}$ are dual orthonormal bases of $\h^*$
and $\h$. Therefore, every lowest weight $H_{\sf k}(W)$-module $M$
is also an $\mathfrak{sl}_2$-module from the category $\mathcal
O$. The element $\mathbf h$ is $W$-invariant and satisfies the
relations

$$[\mathbf h, x]=x, \quad [\mathbf h, y]=-y$$
for all $x \in \h^*$ and $y \in \h$. The spectrum of $\mathbf h$
on $M$ is $\mathbb Z_{\geqslant 0}+h_0$, where $h_0$ is the lowest
weight of $\mathbf h$ on $M$.

\begin{remark} \label{grading}
The generalized eigenspace of $\mathbf h$ with eigenvalue $h_0+i$
consists of all the elements of degree $i$ in $M$.
\end{remark}

As in the case of Lie groups, one can define the category $\O$ of
$H_{\sf k}(W)$-modules.

\begin{definition}{\it The category $\mathcal O$} is the category
of $H_{\sf k}(W)$-modules $M$, such that $M$ is a direct sums of
finite dimensional generalized eigenspaces of $\mathbf h$, and the
spectrum of $\mathbf h$ is bounded from below.
\end{definition}

Standard modules $M_{\sf k} (V)$ and their quotients are objects
of the category $\O$. It is known that $\{L_{\sf k}(V)\}_{V \in
Irrep(W)}$ are precisely the simple objects of $\O$.

\medskip

For $M \in \O$, let $M=\oplus_a M_a$ be its $\mathbf h$-weight
decomposition.

\begin{definition}{\it The character} of $M$ is a formal power
series in $t$ given by
$$
\chi_M(g,t)=\sum_a \Tr \bigl |_{M_a}(g) t^a,
$$
where $g \in W$, and $\Tr \bigl |_{M_a}(g)$ is the trace of
$g$-action on $M_a$.
\end{definition}

The character of the standard module $M_{\sf k}(V)$ is given by
$$
\chi_{M_{\sf k}(V)}(g,t)=\frac{\chi_{V}(g)t^{h_0(V)}}{\det \bigl
|_{\h} (1-gt)},
$$
where $  h_0(V)= \frac{1}{2}\dim \h - \left(\sum_{s \in S} {\sf
k}(s) s\right)\bigr|_{V}$ is the lowest weight of $\mathbf h$ on
$M_{\sf k} (V)$ (see \cite{BEG}, formula (1.5)).

\begin{remark} \label{symmetry}
For any character $\eps: W \to \CC^*$ there is a symmetry between
$H_{\sf k}(W)$ and $H_{\eps {\sf k}}(W)$. Namely, there exists an
isomorphism $H_{\sf k}(W) \stackrel{\sim}{\to} H_{\eps {\sf
k}}(W)$, which acts identically on $\h$ and $\h^*$, and sends $w
\in W$ to $\eps(w)w$. This isomorphism induces an equivalence of
categories of representations of these algebras. Under this
equivalence, $M_{\sf k}(V)$ goes to $M_{\eps {\sf k}}(\eps \ot
V)$.
\end{remark}

\section{Results}

The group $I_2(d)$ acts via reflections on $\h_{\mathbb R} \simeq
\CC$. The subset $S \subset W$ of reflections consists of elements
$s_j: z \mapsto \omega^j\bar{z}$, where $\omega=e^{2\pi i /d}$ and
$1 \leqslant j \leqslant d$.

\subsection{The case $d=2m+1$.}
All the reflections in $S$ are conjugate, and an
$I_2(m)$-invariant function ${\sf k}: S \to \CC$ is a constant
function
 ${\sf k}=c \in \CC$. Denote by $H_c$ the corresponding rational
Cherednik algebra $H_{\sf k}(W)$.

The irreducible representations of $I_2(d)$ are: the trivial
representation $\sf triv$, the sign representation $\sf sgn$, and
two-dimensional representations $\tau_l$, where $1 \leqslant l
\leqslant m$. The representations $\tau_l$ are realized on the
span of $\{z^l, \bar{z}^l\}$.

In what follows, we denote the product $dc$ by $r$.

\begin{theorem}\label{odd1dim} The structure of $M_c(\sf triv)$ for
different values of $c$  is described in the table below.
\begin{center}
{\renewcommand{\arraystretch}{2}
\begin{tabular}{|p{3.2cm}|p{4.9cm}|p{3.5cm}|}
\hline
\centerline{\bf Values of $c$} &
   \bf Structure of the standard\newline module &
   \bf Character of the\newline irreducible quotient \\
\hline
$0<r \equiv \pm l\ (\!\!\!\!\mod d)$,\vspace{5pt}\newline
      where $1 \leqslant l \leqslant m$ &
   $M_c(\sf triv)$ has a submodule isomorphic to $L_c(\tau_l)$.
      The quotient is irreducible of dimension $r^2$ . &
   \centerline{\raisebox{-10pt}{$\displaystyle t^{1-r}
      \frac{\det\bigl|_{\tau_l}(1-gt^{r})}{\det\bigl|_{\h}(1-gt)}$}}\\
\hline
$0<c=n+\frac{1}{2}$,\vspace{5pt}\newline  where $n \in \mathbb Z$
&
   $M_c({\sf triv})$ has a submodule isomorphic to $M_c{(\sf sgn)}$.
      The quotient is irreducible and infinite dimensional.  &
   \centerline{\raisebox{-15pt}{$\displaystyle t^{1-r}
      \frac{1-{\sf sgn}(g)t^{2r}}{\det\bigl|_{\h}(1-gt)}$}}\\
\hline
Other values of $c$ &
   $M_c{(\sf triv})$ is irreducible. &
   \centerline{$\displaystyle t^{1-r}\frac{1}{\det\bigl|_{\h} (1-gt)}$}\\
\hline
\end{tabular}
}
\end{center}
\end{theorem}

\smallskip

The description of the $M_c({\sf sgn})$ follows from theorem
\ref{odd1dim} by symmetry (Remark \ref{symmetry}).

\begin{theorem}\label{odd2dim}
Fix $1 \leqslant l \leqslant m$. The structure of $M_c(\tau_l)$
for different values of $c$ is described in the following table.
\begin{center}
{\renewcommand{\arraystretch}{2}
\begin{tabular}{|p{3.2cm}|p{4.9cm}|p{3.5cm}|}
\hline
\centerline{\bf Values of $c$} &
   \bf Structure of the standard\newline module &
   \bf Character of the\newline irreducible quotient \\
\hline
$0<r \equiv \pm l\ (\!\!\!\!\mod d)$ &
   $M_c(\tau_l)$ has a submodule isomorphic to $M_c{(\sf sgn)}$.
      The quotient is irreducible and infinite dimensional. &
   \centerline{\raisebox{-10pt}{$\displaystyle
      t\frac{\chi_{\tau_l}(g)-t^r}{\det\bigl|_{\h}(1-gt)} {\sf
sgn}(g)$}}\\
\hline
$0>r \equiv \pm l\ (\!\!\!\!\mod d)$ &
   $M_c(\tau_l)$ has a submodule isomorphic to $M_c{(\sf triv)}$.
      The quotient is irreducible and infinite dimensional. &
   \centerline{\raisebox{-10pt}{$\displaystyle t
      \frac{\chi_{\tau_l}(g)-t^{-r}}{\det\bigl|_{\h}(1-gt)}$}}\\
\hline
Other values of $c$ &
   $M_c(\tau_l)$ is irreducible. &
   \centerline{$\displaystyle
      t\frac{\chi_{\tau_j}(g)}{\det\bigl|_{\h} (1-gt)}$}\\
\hline
\end{tabular}
}
\end{center}
\end{theorem}

\smallskip

\begin{remark}
All the standard modules over $H_c$ are of length either $1$ or
$2$.
\end{remark}

\bigskip

\subsection{The case  $d=2m$.}
There are two conjugacy classes of reflections:  $\{s_{2i}\}_{1
\leqslant i \leqslant m}$ and $\{ s_{2i-1}\}_{1 \leqslant i
\leqslant m }$. Hence a $W$-invariant function ${\sf k}:S \to
\CC$ is a pair $(k_1, k_2)$, where $k_1={\sf k}(s_{2i+1})$ and
$k_2 = {\sf k}(s_{2i})$. Let $H_{k_1, k_2}$ denote the
corresponding rational Cherednik algebra $H_{\sf k} (W)$.

The group $I_2(d)$ has four irreducible one-dimensional
representations: ${\sf triv},$ ${\sf sgn},$ $\eps_1,$
$\eps_2=\eps_1 \ot {\sf sgn}$, and $(m-1)$ two-dimensional
representations $\tau_l$, $1 \leqslant l < m$. The representation
$\eps_1$ is defined by
$$
\eps_1(s_{2i+1})=1, \quad \eps_1(s_{2i})=-1,
$$
and $\tau_l$'s are realized on $\linob\{z^l, \bar{z}^l\}$.

 \begin{remark}
 Note that ${\sf sgn} \ot \tau_l =\tau_l$ and $\eps_i \ot \tau_l
 =\tau_{m-l}$ for $i=1,2$.
 \end{remark}

\def\wt#1{\widetilde{#1}}
The structure of the standard $H_{k_1, k_2}$-modules depends on
the location of the point $(k_1,k_2)$ with respect to certain
exceptional lines on the complex plane of parameters $(k_1,k_2)$.
We distinguish the four families $E^+_r$, $E^-_r$, $L^1_i$,
$L^2_i$ of parallel exceptional lines: \ben
\item[$\bullet$] for every integer $r$, non divisible by $m$, let
$$E^+_r := \{(k_1,k_2) : k_1+k_2 = r/m\}, \qquad\qquad
  E^-_r := \{(k_1,k_2) : k_1-k_2 = r/m\}\ ;
$$
\item[$\bullet$] for every integer $i$, let
$$L^1_i := \{(k_1,k_2) : k_1 = i+\half\}, \qquad\qquad
  L^2_i := \{(k_1,k_2) : k_2 = i+\half\}\ .
$$
\een

\medskip
By symmetry reasons, in order to describe the standard modules
induced from one-dimensional representations, it is enough to
describe $M_{k_1,k_2}({\sf triv})$.

\def\ChTriv{\chi_{\raisebox{-5pt}{\!\!
     $\scriptstyle L_{k_1,k_2}({\sf triv})$}}\!(g,t)}
\def\ChTaul{\chi_{\raisebox{-5pt}{\!\!
     $\scriptstyle L_{k_1,k_2}(\tau_l)$}}\!(g,t)}
\begin{theorem}\label{even1dim}
Depending on the location of the point $(k_1,k_2)$, there are
following seven cases for the  standard module $M_{k_1,k_2}({\sf
triv})$. \ben
\item[({\bf i})] {\bf $(k_1,k_2) \in E_r^+$ for some $r>0$, and
           $(k_1,k_2) \notin L_i^1, L^2_i$ for all $i>0$.}

In this case determine an integer $l$ in the range $1\leqslant l <m$
from the congruence
$r\equiv \pm l\ (\!\!\!\!\mod 2m)$. Then $M_{k_1,k_2}({\sf triv})$
contains a submodule isomorphic to $L_{k_1,k_2}(\tau_l)$. The
quotient is irreducible and has the character
$$\ChTriv \quad =\quad t^{1-r}
   \frac{\det\bigl|_{\tau_l}(1-gt^r)}{ \det\bigl|_{\h} (1-gt)}\ .
$$

In particular, its dimension is $r^2$.

\medskip
\item[({\bf ii})] {\bf $(k_1, k_2) \in L_i^1$ for some $i>0$, and
  $(k_1,k_2) \notin E_r^+, L^2_{i'}$ for  all $r, i'>0$.}

In this case the module
$M_{k_1,k_2}({\sf triv})$ contains a submodule isomorphic to
$L_{k_1,k_2}(\eps_2) \simeq M_{k_1,k_2}(\eps_2)$. The quotient is
irreducible and infinite dimensional. The character is given by the
formula
$$\ChTriv \quad =\quad t^{1-m(k_1+k_2)}
   \frac{1-t^{2mk_1}\chi_{\eps_2}(g)}{\det\bigl|_{\h}(1-gt)}\ .
$$

\medskip
\item[({\bf iii})] {\bf $(k_1, k_2) \in  L^2_{i'}$ for $i'>0$, and
  $(k_1,k_2)\notin E_r^+, L_i^1$ for all $r,i >0$.}

In this case the module
$M_{k_1,k_2}({\sf triv})$ contains a submodule isomorphic to
$L_{k_1,k_2}(\eps_1) \simeq M_{k_1,k_2}(\eps_1)$. The quotient is
irreducible and infinite dimensional. The character is given by the
formula
$$\ChTriv \quad =\quad t^{1-m(k_1+k_2)}
   \frac{1-t^{2mk_2}\chi_{\eps_1}(g)}{\det\bigl|_{\h}(1-gt)}\ .
$$

\medskip
\item[({\bf iv})] {\bf $(k_1,k_2) \in E_r^+ \cap L_i^1$
   (resp. $(k_1,k_2) \in E_r^+ \cap L^2_{i'}$) for some $r,i,i'>0$, and
   $k_1< k_2$ (resp. $k_2 < k_1$).}

Determine an integer $l$ from the same conditions
$r\equiv \pm l\ (\!\!\!\!\mod 2m)$ and $1\leqslant l <m$.
We have the following inclusions of the submodules
$$
M_{k_1, k_2}({\sf triv}) \supset M_{k_1,k_2}(\eps_2)
  \supset L_{k_1,k_2}(\tau_l)
$$
$$\mbox{{\rm (}resp.}\quad M_{k_1,k_2}({\sf triv}) \supset
   M_{k_1,k_2}(\eps_1) \supset L_{k_1,k_2}(\tau_l)\ \mbox{\rm )}\ .
$$
The quotient $M_{k_1,k_2}(\eps_2) / L_{k_1,k_2}(\tau_l)$\quad {\rm
(}resp. $M_{k_1,k_2}(\eps_1) / L_{k_1,k_2}(\tau_l)${\rm )}, is
irreducible. The  quotient $Q_{k_1,k_2}=M_{k_1, k_2}({\sf triv}) /
M_{k_1,k_2}(\eps_2)$ {\rm (}resp. $Q_{k_1,k_2}=M_{k_1, k_2}({\sf
triv}) / M_{k_1,k_2}(\eps_1)${\rm )} contains a submodule
$P_{k_1,k_2}$ isomorphic to \linebreak $L_{k_1,k_2}(\eps_1)$ {\rm
(}resp. $L_{k_1,k_2}(\eps_2)${\rm )}. However, the double quotient
$Q_{k_1,k_2}/P_{k_1,k_2}$ is irreducible, i.e. isomorphic to
$L_{k_1,k_2}({\sf triv})$.

In both cases the character of $L_{k_1,k_2}({\sf triv})$ is given by
the formula
$$\qquad\qquad\ChTriv \quad =\quad t^{1-r}
   \frac{1+t^{2r}{\sf sgn}(g)-t^{2mk_1}\eps_2(g)-t^{2mk_2}\eps_1(g)}
   {\det\bigl|_{\h}(1-gt)}\ .
$$
In particular, its dimension is equal to $4m^2k_1k_2$.

\medskip
\item[({\bf v})] {\bf $(k_1,k_2)\in E_r^+ \cap L_i^1$
   (resp. $(k_1,k_2) \in E_r^+ \cap L^2_{i'}$) for some $r,i,i'>0$, and
   $k_1 > k_2$ (resp. $k_2 > k_1$).}

Again, let $l$ be an integer, satisfying the conditions $r\equiv
\pm l\ (\!\!\!\!\mod 2m)$ and $1\leqslant l <m$. Then $M_{k_1,
k_2}({\sf triv})$ contains a submodule isomorphic to
$L_{k_1,k_2}(\eps_2)$ {\rm (}resp. $L_{k_1,k_2}(\eps_1)${\rm )}.
The quotient $Q_{k_1,k_2}=M_{k_1, k_2}({\sf triv}) /
L_{k_1,k_2}(\eps_2)$ {\rm (}resp. $Q_{k_1,k_2}=M_{k_1, k_2}({\sf
triv}) / L_{k_1,k_2}(\eps_1)${\rm )} contains a submodule
$P_{k_1,k_2}$ isomorphic to $L_{k_1,k_2}(\tau_l)$. The quotient
$Q_{k_1,k_2}/P_{k_1,k_2}$ is irreducible and has the character
$$\ChTriv \quad =\quad t^{1-r}
   \frac{\det\bigl|_{\tau_l}(1-gt^r)}{ \det\bigl|_{\h} (1-gt)}\ .
$$
In particular, its dimension is $r^2$.

\medskip
\item[({\bf vi})] {\bf $(k_1,k_2) \in L_i^1 \cap L^2_{i'}$ for
   some $i,i'>0$. }

In this case the module $M_{k_1,k_2}({\sf triv})$ contains a
submodule isomorphic to $L_{k_1,k_2}({\sf sgn}) \simeq
M_{k_1,k_2}({\sf sgn})$. The quotient 
$Q_{k_1,k_2}=M_{k_1,k_2}({\sf triv}) / L_{k_1,k_2}({\sf sgn})$
contains two submodules $P^1_{k_1,k_2} \simeq L_{k_1,k_2}(\eps_1)$
and $P^2_{k_1,k_2} \simeq L_{k_1,k_2}(\eps_2)$, such that
$P^1_{k_1,k_2} \cap P^2_{k_1,k_2} =0$. The quotient
$Q_{k_1,k_2}/(P^1_{k_1,k_2} \oplus P^2_{k_1,k_2})$ is irreducible,
i.e. isomorphic to $L_{k_1,k_2}({\sf triv})$ with the character
$$\ChTriv \quad =\quad t^{1-r}
   \frac{1+t^{2r}{\sf sgn}(g)-t^{2mk_1}\eps_2(g)-t^{2mk_2}\eps_1(g)}
        {\det\bigl|_{\h}(1-gt)}\ .
$$
In particular, it has dimension $ 4 m^2 k_1k_2$.

\medskip
\item[({\bf vii})] {\bf For all other values of $(k_1,k_2)$}

the  module $M_{k_1,k_2}({\sf triv})$ is irreducible and has the
character
$$\chi_{\raisebox{-5pt}{\!\!
     $\scriptstyle L_{k_1,k_2}({\sf triv})$}}\!(g,t)\quad =\quad
     t^{1-r}\frac{1}{\det\bigl|_{\h}(1-gt)}\ .
$$
\een
\end{theorem}

\bigskip
The structure of the standard modules, induced from the
two-dimensional representations $\tau_l$ of the group $W=I_2(d)$
is described by the following theorem.

\bigskip
\begin{theorem}\label{even2dim}
Fix $1 \leqslant l < m$ and consider the two-dimensional
representation $\tau_l$ of $I_2(d)$. Then, depending on the
location of the point $(k_1,k_2)$, there are following four cases
for the  standard module $M_{k_1,k_2}(\tau_l)$. \ben
\item[({\bf i})] {\bf $(k_1, k_2) \in E_r^+$ for some
      $r \equiv \pm l\ (\!\!\!\!\mod 2m)$, and
      $(k_1, k_2) \notin E^-_{r'}$ for all
      $r' \equiv m \pm l\ (\!\!\!\!\mod 2m)$.}

The module
$M_{k_1, k_2} (\tau_l)$ contains a submodule $I_{k_1, k_2}$
generated by one vector of degree $|r|$.

If $r>0$, then $I_{k_1,k_2} \simeq M_{k_1, k_2}(\sf sgn)$.
If $r<0$, then $I_{k_1, k_2} \simeq M_{k_1, k_2}(\sf triv)$.

The quotient $M_{k_1,k_2}(\tau_l)/I_{k_1, k_2}$ is irreducible and
infinite
dimensional. If $r>0$, its character is equal to
$$\ChTaul \quad =\quad t
   \frac{\chi_{\tau_l}(g)-t^r}{\det\bigl|_{\h}(1-gt)}{\sf
sgn}(g)\ .
$$
For $r<0$, the character can be easily obtained using symmetry
arguments.

\medskip
\item[({\bf ii})] {\bf $(k_1, k_2) \notin E_{r}$ for all
      $r \equiv \pm l\ (\!\!\!\!\mod 2m)$, but
      $(k_1, k_2) \in E^-_{r'}$ for some
      $r' \equiv m \pm l\ (\!\!\!\!\mod 2m)$.}

The module $M_{k_1, k_2}(\tau_l)$ contains a
submodule $\wt{I}_{k_1, k_2}$ generated by one vector
of degree $|r'|$.

If $r'>0$, then $\wt{I}_{k_1,k_2} \simeq M_{k_1,k_2}(\eps_2)$.
If $r'<0$, then $\wt{I}_{k_1, k_2}\simeq M_{k_1, k_2}(\eps_1)$.

The quotient $M_{k_1, k_2}(\tau_l)/\wt{I}_{k_1, k_2}$ is
irreducible. If $r'>0$, its character is equal to
$$\ChTaul \quad =\quad t
   \frac{\chi_{\tau_l}(g)-t^{r'}}{ \det\bigl|_{\h}(1-gt)}\eps_2(g)\
.
$$
For $r'<0$, the character can be easily obtained using symmetry
arguments.

\medskip
\item[({\bf iii})] {\bf $(k_1,k_2) \in E_r^+ \cap E^-_{r'}$, for some
      $r \equiv \pm l\ (\!\!\!\!\mod 2m)$ and
      $r' \equiv m \pm l\ (\!\!\!\!\mod 2m)$.}

The module $M_{k_1, k_2} (\tau_l)$ contains two submodules:
$I_{k_1, k_2}$ and $\wt{I}_{k_1, k_2}$.

If\quad $\min\{|k_1|, |k_2|\} \in \mathbb{Z}+\half$,\quad then one
of these submodules contains another one; otherwise they intersect
by $0$. The quotient of $M_{k_1, k_2} (\tau_l)$ by those two
submodules is irreducible in both cases. If $r,r'>0$, then the
character is
$$\ChTaul \quad =\quad t
   \frac{\chi_{\tau_l}(g)-t^r}{\det\bigl|_{\h}(1-gt)}{\sf sgn}(g)
$$
in the first case, and
$$\ChTaul \quad =\quad t
   \frac{\chi_{\tau_l}(g)-t^r{\sf sgn}(g) -t^{r'}\eps_2(g)}
        {\det\bigl|_{\h}(1-gt)}
$$
in the second case.\\
For other combinations of signs of $r$ and $r'$, the characters
can be easily obtained using symmetry arguments.

\medskip
\item[({\bf iv})] {\bf For all other values of $(k_1,k_2)$}

the  module $M_{k_1,k_2}(\tau_l)$ is irreducible and has the character
$$\chi_{\raisebox{-5pt}{\!\!
     $\scriptstyle L_{k_1,k_2}(\tau_l)$}}\!(g,t)\quad =\quad
   t\frac{\tau_l(g)}{\det\bigl|_{\h}(1-gt)}\ .
$$
\een
\end{theorem}


\section{Examples.}
In order to illustrate the Theorems \ref{even1dim} and
\ref{even2dim}, let us look at $W=I_2(8)$. The Theorem
\ref{even1dim} gives us the following picture for
$M_{k_1,k_2}({\sf triv})$ in $(k_1,k_2)$-plane.

\mbox{
\begin{picture}(140,150)(0,0)
\put(2,2){\epsfxsize=130pt \epsfbox{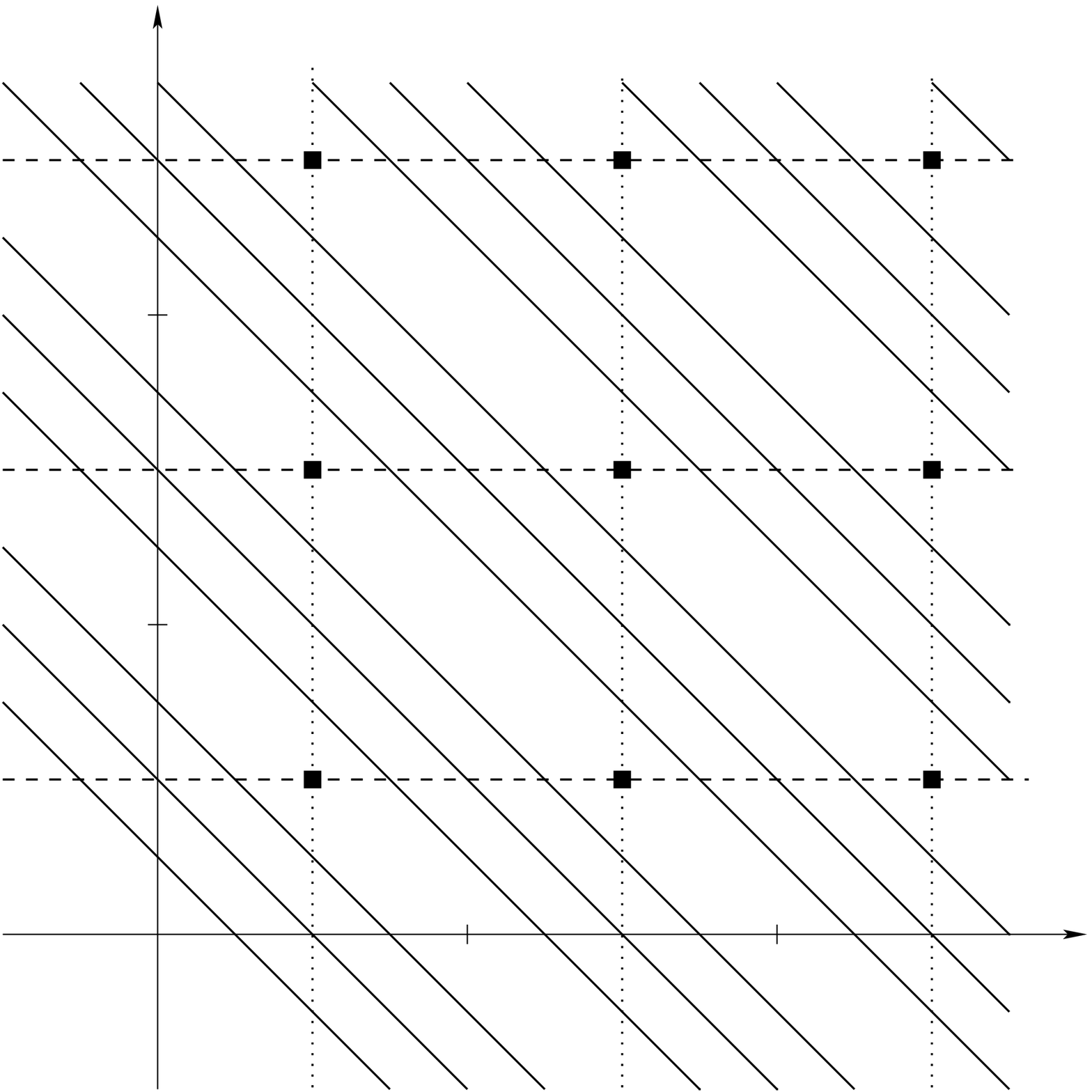}}
 \put(130,14){\mbox{$_{k_1}$}}
 \put(10, 133){\mbox{$_{k_2}$}}
 \put(56,14){\mbox{$_1$}}
 \put(14,57){\mbox{$_1$}}
 \put(93,14){\mbox{$_2$}}
 \put(14,94){\mbox{$_2$}}
\end{picture}
} \mbox{
\begin{picture}(140,150)(0,0)
\put(2,20){\epsfxsize=15pt \epsfbox{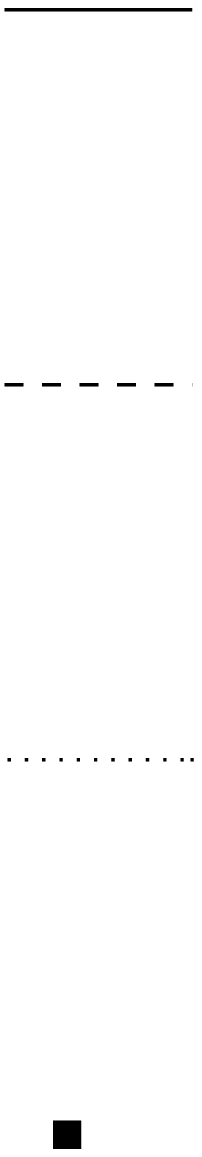}}
 \put(60,130){\mbox{Notations:}}
 \put(22,105){\mbox{-- $M_{k_1,k_2}({\sf triv}) \supset L_{k_1,k_2} (\tau_i)$ for some $i$},}
 \put(5,92){\mbox{the quotient is finite dimensional;}}
 \put(22,76){\mbox{-- $M_{k_1,k_2}({\sf triv}) \supset M_{k_1,k_2} (\eps_1)$, quotient}}
 \put(5,63){\mbox{is generically infinite dimensional;}}
 \put(22,47){\mbox{-- $M_{k_1,k_2}({\sf triv}) \supset M_{k_1,k_2} (\eps_2)$, quotient}}
 \put(5,34){\mbox{is generically infinite dimensional;}}
 \put(22,19){\mbox{-- $M_{k_1,k_2}({\sf triv})$ contains $M_{k_1,k_2} (\eps_1)$ and}}
 \put(5,6){\mbox{$M_{k_1,k_2}(\eps_2)$, quotient is finite dimensional.}}
\end{picture}
}

For all the points outside the indicated lines, $M_{k_1,k_2}({\sf
triv})$ is irreducible.

Using this information and the symmetry arguments, one can find
all the values of $(k_1,k_2)$, for which there exist finite
dimensional representation of $H_{k_1,k_2}$:
$$
\mbox{
\begin{picture}(255,220)(0,0)
\put(2,2){\epsfxsize=250pt \epsfbox{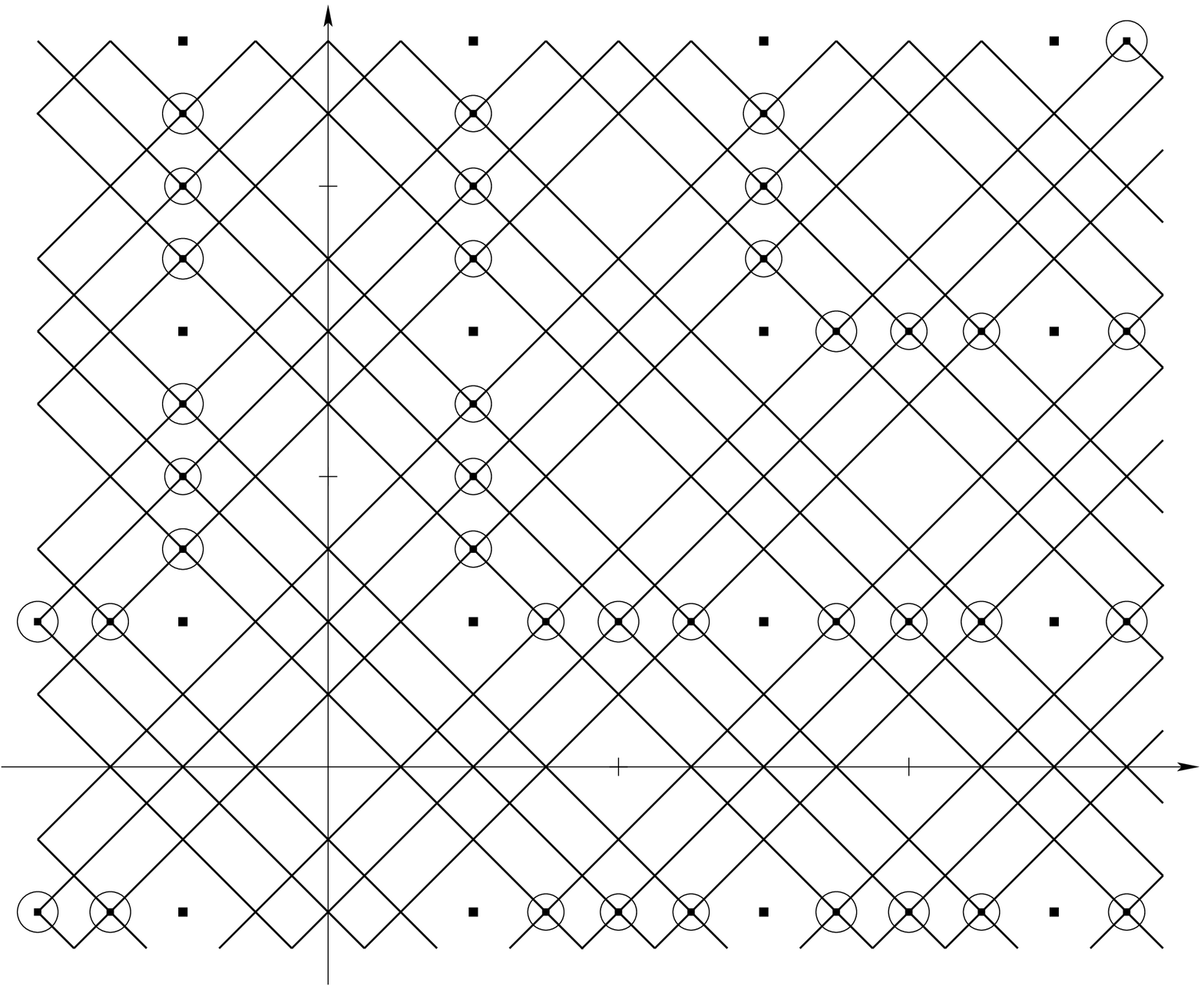}}
 \put(250,40){\mbox{$_{k_1}$}}
 \put(60, 205){\mbox{$_{k_2}$}}
 \put(129,40){\mbox{$_1$}}
 \put(64,108){\mbox{$_1$}}
 \put(189,40){\mbox{$_2$}}
 \put(64,168){\mbox{$_2$}}
\end{picture}
}
$$

The following table lists the finite dimensional irreducible
representations for $k_1,k_2 >0$.

$$
{\renewcommand{\arraystretch}{2}%
\begin{tabular}{|p{5cm}|c|p{3.5cm}|c|}
\hline
\multicolumn{2}{|c|}{\bf Values of parameters}
   & {\bf Finite dimensional\newline irreducible modules}
   & {\bf Dimension} \\
\hline
 \multicolumn{2}{|l|}{$(k_1,k_2)$ is generic point on $E_r^+$}  &
   \hspace{1cm}$L_{k_1,k_2}({\sf triv})$& $r^2$ \\
\hline \makebox(137,0){\raisebox{-23pt}{$(k_1,k_2)$ is generic
point on
                                 $E^-_{r'}$}} &
      $r'>0$ & \hspace{1cm}$L_{k_1,k_2}(\eps_1)$&
      \makebox(33,0){\raisebox{-23pt}{${r'}^2$}} \\
 \hhline{|~-|-|~|}
 & $r'<0$ &
 \hspace{1cm}$L_{k_1,k_2}(\eps_2)$&  \\
 \hline
\makebox(125,0){\raisebox{-28pt}{\parbox{125pt}{
       $(k_1,k_2) \in E_r^+ \cap E_{r'}^-$, \vspace{3pt}
        denoted by $\odot$ on the picture above }}}&
  $r'>0$ & $L_{k_1,k_2 }({\sf triv})$, $L_{k_1,k_1}(\eps_1)$ &
      \makebox(33,0){\raisebox{-23pt}{$r^2-{r'}^2$, ${r'}^2$}} \\
 \hhline{|~-|-|~|}
 & $r'<0$ & $L_{k_1,k_2 }({\sf triv})$, $L_{k_1,k_1}(\eps_2)$ & \\
\hline \makebox(137,0){\raisebox{-28pt}{\parbox{137pt}{
       All the other $(k_1,k_2) \in E_r^+ \cap E_{r'}^-$ }}}&
  $r'>0$ & $L_{k_1,k_2 }({\sf triv})$, $L_{k_1,k_1}(\eps_1)$ &
      \makebox(33,0){\raisebox{-23pt}{$r^2$, ${r'}^2$}} \\
 \hhline{|~-|-|~|}
 & $r'<0$ & $L_{k_1,k_2 }({\sf triv})$, $L_{k_1,k_1}(\eps_2)$ & \\
\hline
 \multicolumn{2}{|l|}{Isolated points on the picture above}  &
              \hspace{1cm}$L_{k_1,k_2}({\sf triv})$& $ 64 k_1k_2$\\
\hline
\end{tabular}}
$$

\begin{remark}
If $(k_1,k_2)$ lies on one of the exceptional lines, the category
$\mathcal O$ is not semisimple.

If $(k_1,k_2)$ is one of the points denoted by $\odot$ in the
picture above, the category of finite dimensional
$H_{k_1,k_2}$-modules is not semisimple. For example, in the case
of $k_1 > k_2 >0 \quad \text{Ext}^1(L_{k_1,k_2}({\sf triv}),
L_{k_1,k_2}(\eps_1)) \ne 0$.
\end{remark}

Applying the theorem \ref{even2dim} to  $M_{k_1, k_2} (\tau_1)$,
we will get the following picture:

\mbox{
\begin{picture}(210,210)(0,0)
\put(2,2){\epsfxsize=200pt \epsfbox{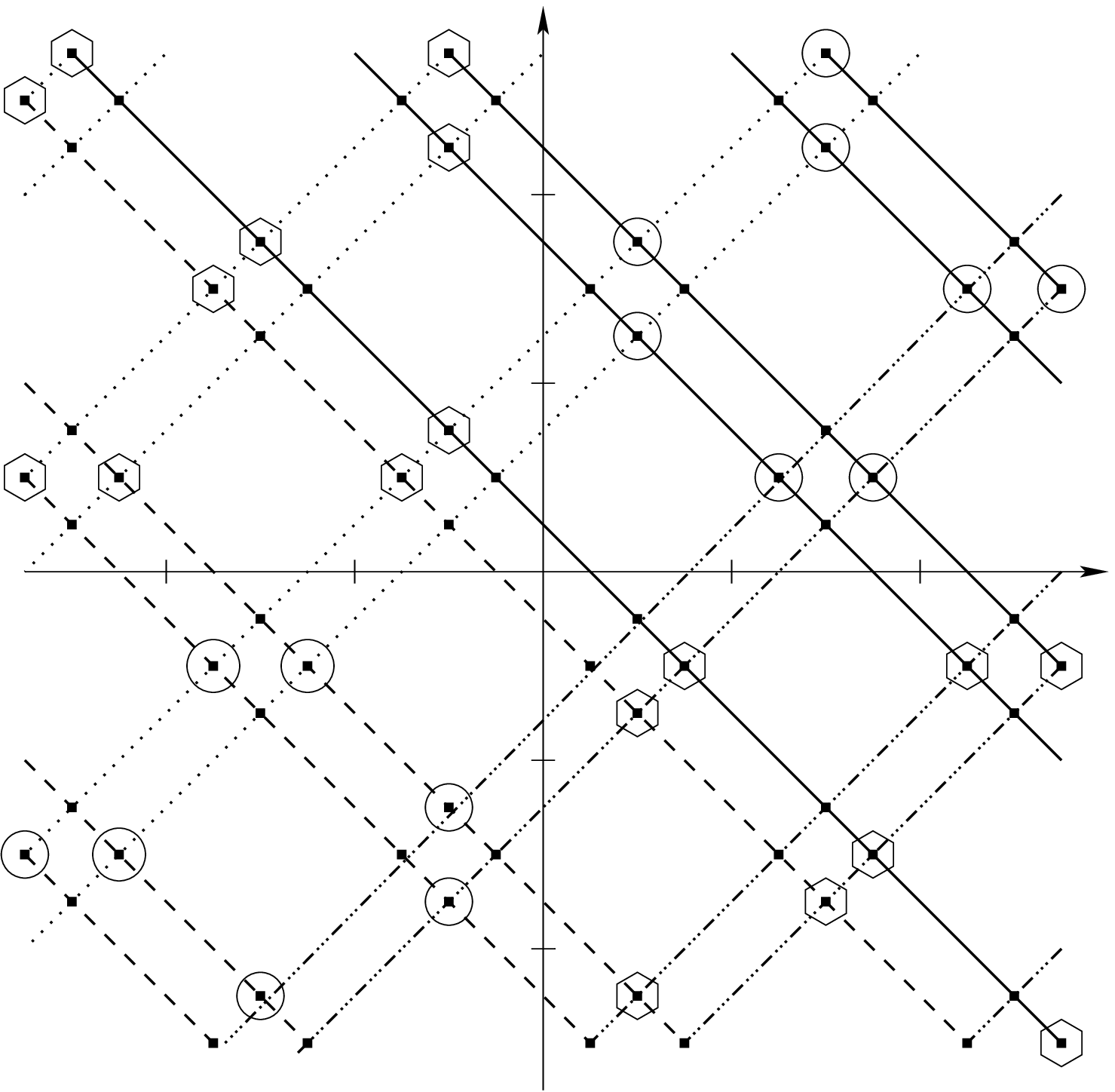}}
 \put(200,90){\mbox{$_{k_1}$}}
 \put(90, 200){\mbox{$_{k_2}$}}
 \put(130,90){\mbox{$_1$}}
 \put(93,125){\mbox{$_1$}}
\end{picture}
}
 \mbox{
\begin{picture}(100,210)(0,0)
\put(2,2){\epsfxsize=25pt \epsfbox{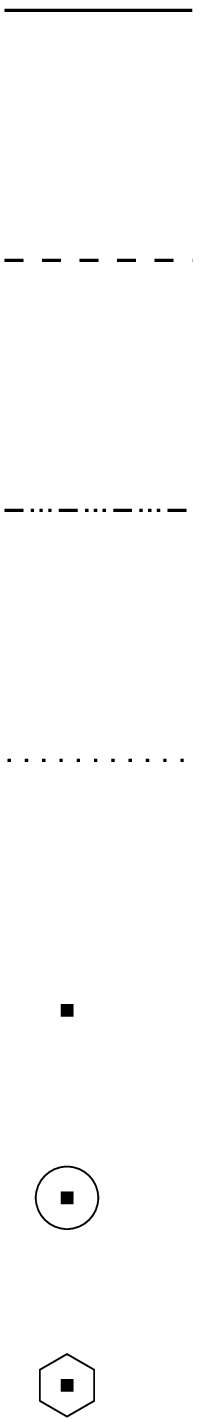}}
 \put(40,200){\mbox{ Notations:}}
 \put(30,177){\mbox{-- $M_{k_1,k_2}(\tau_1) \supset I_{k_1,k_2}$}}
 \put(10,164){\mbox{and $I_{k_1,k_2}\simeq M_{k_1,k_2}({\sf sgn})$};}
 \put(30,146){\mbox{-- $M_{k_1,k_2}(\tau_1) \supset I_{k_1,k_2}$ }}
 \put(10,133){\mbox{and $I_{k_1,k_2}\simeq M_{k_1,k_2}({\sf triv})$};}
 \put(30,114){\mbox{-- $M_{k_1,k_2}(\tau_1) \supset \widetilde{I}_{k_1,k_2}$ }}
 \put(10,101){\mbox{and $\widetilde{I}_{k_1,k_2}\simeq M_{k_1,k_2}({\eps_2})$};}
 \put(30,83){\mbox{-- $M_{k_1,k_2}(\tau_1) \supset \widetilde{I}_{k_1,k_2}$ }}
 \put(10,70){\mbox{and $\widetilde{I}_{k_1,k_2}\simeq M_{k_1,k_2}({\eps_1})$};}
 \put(18,51){\mbox{-- $I_{k_1,k_2} \cap \widetilde{I}_{k_1,k_2}=0$;}}
 \put(18,28){\mbox{-- $I_{k_1,k_2} \subset \widetilde{I}_{k_1,k_2}$;}}
 \put(18,3){\mbox{-- $I_{k_1,k_2} \supset \widetilde{I}_{k_1,k_2}$.}}
\end{picture}
}

As before, for all the points outside the indicated lines,
$M_{k_1, k_2} (\tau_1)$ is irreducible.

\section{Preparations for the proofs.}

In this section we collect known results, which we need for the
proofs.

\subsection{The Knizhnik-Zamolodchikov functor.}\label{KZ}
We will briefly remind the construction of the $\sf KZ$-functor
which, maps certain modules over the Cherednik algebra to modules
over the Artin braid group. A more detailed discussion on this
subject can be found in \cite{BEG1}, \cite{GGOR} and \cite{BMR}.

\subsubsection{Hecke algebras and the Knizhnik-Zamolodchikov functor.}
Let us first recall the definition of Hecke algebras associated to
a reflection group $W$.

Let $\h_{reg}=\h \backslash  \bigcup_{s \in S} H_s $. The
$W$-action on $\h$ preserves $\h_{reg}$. Fix $x_0 \in \h_{reg}$.

\begin{definition}
The {\it Artin braid group} associated to $W$ is
$B_W=\pi_1(\h_{reg}/W, x_0)$.
\end{definition}

One can describe $B_W$ using generators and relations as follows
(see \cite{Br}). Fix a chamber $C$ in $\h$, and let $\Sigma$ be
the set of reflections with respect to the walls of this chamber.
The generators $T_s$ are numbered by element $s \in \Sigma$. An
element $T_s$ corresponds to the monodromy around the wall of $C$,
fixed by $s$. The relations are
$$\underbrace{T_s T_t T_s \dots}_{\mbox{\scriptsize{$m(s,t)$ factors}}}
 = \underbrace{T_t T_s T_t \dots}_{\mbox{\scriptsize{$m(s,t)$ factors}}}
 \text{ for all } s, t \in S,$$
where $m(s,t)$ is the order of the element $st$.

Let $\{q_s\}_{s \in \Sigma}$ be the set of parameters such that
$q_s=q_t$ if $s$ and $t$ are conjugated to each other.
\begin{definition}
The {\it Hecke algebra} $\mathcal H (\{q_s\})$ of a reflection
group $W$ with parameters $q_s$ is the quotient of the group
algebra of the Artin braid group $B_W$ by the relations
$$
(T_s-1)(T_s+q_s)=0.
$$
\end{definition}

Now let us define the Knizhnik-Zamolodchikov functor $\sf KZ$. For
any $H_{\sf k}$-module $M$, finitely generated over $\mathbb C
[\h]$, define its localization $M_{loc}$ by
$$M_{loc}=M
\otimes_{\mathbb C[\h]} \mathbb C[\h_{reg}].$$
 This module has a natural action of ${H_{\sf k \hspace{2pt}}}_{loc}$.

The polynomial representation  (see section \ref{standmod}) gives
an injective homomorphism $H_{\sf k} \hookrightarrow W \ltimes
 \mathcal D(\h_{reg})$, which extends to an isomorphism
${H_{\sf k \hspace{2pt}}}_{loc} \stackrel{\sim}{\rightarrow} W
\ltimes \mathcal D(\h_{reg})$. Hence the module $M_{loc}$ has a
structure of $W \ltimes \mathcal D(\h_{reg})$-module, and
therefore, a structure of $W$-equivariant
 $\mathcal D$-module on $\h_{reg}$.

Since $M$ is finitely generated over $\CC [\h]$, the module
$M_{loc}$ is a vector bundle over $\h_{reg}$. This vector bundle
is equipped with a connection $\nabla$, coming from the $\mathcal
D$-module structure. The horizontal sections of this connection
define a monodromy representation
$$\rho_M :\pi_1(\h_{reg}/W, x_0) \to GL(\mathbb C^N),$$
where $N$ is the rank of $M_{loc}$.

\begin{definition} A functor ${\sf KZ}:M \mapsto \rho_M$ is called {\it the
Knizhnik-Zamolodchikov functor}.
\end{definition}
\noindent
$$
{\sf KZ}:
 \left(
\begin{array}{c}
H_{\sf k}(W) \text{-modules finitely generated}\\
\text{over } \CC [\h]
\end{array} \right)
 \to
 \left(
\begin{array}{c}
\text{Finite dimensional} \\
B_W \text{-modules}
\end{array} \right)
$$

In case of dihedral groups the images of the standard modules
under the action of ${\sf KZ}$-functor were studied in details by
C.~Dunkl in {\cite{Du1}}.

\subsubsection{Properties of the $\sf
KZ$-functor.}\label{KZproperties}

Let $M$ be an $H_{\sf k}$-module. If $q_s=e^{-2\pi i {\sf k}(s)}$
for all $s$, then the action of $B_W$ on ${\sf KZ}(M)$ factors
through the action of $\mathcal H(\{q_s\})$. If $M$ is the
standard module corresponding to a representation $V$ of the group
$W$, the obtained representation of $\mathcal H (\{q_s\})$ is the
one corresponding to $V$ via Tits' Deformation Theorem.

 The following facts about the $\sf KZ$-functor are important for us.
\begin{proposition}{\rm \cite{GGOR}} The functor $\sf KZ$ is
exact.
\end{proposition}

\begin{theorem}{\rm (\cite{GGOR}, Theorem 5.13)} The functor ${\sf KZ}:\mathcal
O \to B_W$-mod factors through a functor ${\sf KZ}: \mathcal O/
\mathcal O_{tor} \to \mathcal H(\{q_s\})$-mod, where $\mathcal
O_{tor}$ is the category of modules, supported on  the union of
reflection hyperplanes.\label{equiv1}
\end{theorem}

\begin{theorem}{\rm (\cite{BEG1} Lemma 2.10)} Let $N$ be an $H_{\sf k}$-module,
which is torsion free over $\CC [\h]$. Then for any $M \in
\mathcal O$ the canonical map
$$
\Hom_{H_{\sf k}}(M,N) \rightarrow \Hom_{\mathcal H} ({\sf KZ}(M),
{\sf KZ}(N))
$$
is injective. \label{equiv2}
\end{theorem}

\subsection{Singular vectors in the polynomial representation.}
Let $M$ be a module from the category $\mathcal O$. Denote by $U$
the space of singular vectors in $M$. This space is nonzero, since
the action of $y$'s decreases the eigenvalues of $\mathbf h$ by
$1$, and the spectrum of $\mathbf h$ is bounded from below. From
the defining relations of Cherednik algebra follows that $U$ is
preserved under the $W$-action.

\begin{proposition} \label{irred}
A module $M \in \mathcal O$ is irreducible if and only if $U$
spans $M$ and is an irreducible representation of $W$.
\end{proposition}

\begin{flushright}
$\Box$
\end{flushright}

Hence, in order to describe the structure of standard modules, it
is important to study $U$.

Denote by $U' \subset U$  the space of all singular vectors of
positive degree. For the polynomial representation, this space was
studied by C. Dunkl, E. Opdam and M.F.E. de Jeu.

\begin{theorem}\label{oddsingvect}{\rm (see \cite{DJO})}
Consider the Cherednik algebra $H_c(I_2(d))$, where $d=2m+1$.
\begin{enumerate}
\item{Let $c=n \pm l/d >0$, where n and l are integers and $1 \leqslant l \leqslant m$.
Then $U' \subset M_c({\sf triv})$ is isomorphic to $\tau_l$ as a
representation of $W$. This space sits in degree $dn+l$. It is
spanned by
$$
f=\oint \left( \frac{1}{w}\right)^{n+\frac{2l}{d}}z^l
(w-\bar{z}^d)^{n+\frac{l}{d}} (w-z^d)^{n+\frac{l}{d}-1}dw
$$
and its conjugate.}
\item{Let $c=n+1/2>0$. Then $U' \subset M_c({\sf triv})$ is isomorphic to
the sign representation of  $W$, and is spanned by the singular
vector $g=(z^d - \bar{z}^d)^{2n+1}$.}
\item{For all other values of $c$, there are no singular vectors of
positive degree in $M_c({\sf triv})$.}
\end{enumerate}
\end{theorem}

\begin{theorem}\label{evensingvect}{\rm (see \cite{DJO})} Consider
the Cherednik algebra $H_{k_1,k_2}(I_2(2m))$. The space $U'$ is
nonzero only in the following cases.

\ben
\item{If $(k_1, k_2) \in E_r^+$ for some $0<r \equiv \pm l \text{ (mod } 2m$) where $1\leqslant l <m$,
then $U' \subset M_{k_1,k_2}(\sf triv)$ contains a
$W$-subrepresentation isomorphic to $\tau_l$ and consisting of
vectors of degree $r$.}

\item{If $(k_1, k_2) \in L_i^1$ for
some $i \geqslant 0$, then $U' \subset M_{k_1, k_2}(\sf triv)$
contains a $W$-subrepresentation, isomorphic to $\eps_2$ and
spanned by the singular vector $(z^{m}-\bar{z}^{m})^{2i+1}$. }

\item{If $(k_1,k_2) \in L_{i'}^2$ for some $i' \geqslant 0$,
then $U' \subset M_{k_1, k_2}(\sf triv)$ contains a
$W$-subrepresentation isomorphic to $\eps_1$ and spanned by the
singular vector $(z^{m}+\bar{z}^{m})^{2i+1}$.}

\item{If $(k_1,k_2) \in L_{i}^1 \cap L^2_{i'}$ for some $i,
i' \geqslant 0$, then $U' \subset M_{k_1, k_2}(\sf triv)$ is
isomorphic to $\eps_1 \oplus \eps_2 \oplus {\sf sgn}$ as a
$W$-representation. It is spanned by singular vectors
$(z^{m}+\bar{z}^{m})^{2i+1}$, $(z^{m}-\bar{z}^{m})^{2i+1}$, and
$(z^{m}-\bar{z}^{m})^{2i+1}(z^{m}+\bar{z}^{m})^{2i'+1}$.}

\een

Note that if $(k_1,k_2)$ belongs to the intersection of $E_r^+$
with some of the lines $L_i^1$ {\rm (}resp. $L_i^2$ {\rm )}, then
$U'$ is the direct sum of two irreducible representations of $W$,
which are described in the cases {\rm (1)} and {\rm (2)} {\rm
(}resp. {\rm (3))}. The intersection of lines $L_i^1$ and
$L_{i'}^2$ is described separately in the case {\rm (4)}, because
an additional singular vector appears.

\end{theorem}

\subsection{The lowest $\mathbf h$-weights of standard modules.}

Let $V \subset U \subset M$ be an irreducible
$W$-subrepresentation, consisting of singular vectors. Then $V$ is
contained in the generalized eigenspace of $\mathbf h$ with
eigenvalue $h_0(V)$. Indeed, the ideal generated by $V$ is an
$H_{\sf k}(W)$-module of lowest weight $V$. Hence, the smallest
eigenvalue of $\mathbf h$ on this module is $h_0(V)$, and it is
achieved on the generating vectors, i.e. on $V$.

The knowledge of  the spectrum of $\mathbf h$ will help us to
understand $U$ as a representation of $W$.

\begin{proposition}\label{lowesweights}
\ben
\item Let $d=2m+1$. The lowest weights of standard
$H_c$-modules are given in the table below:
\begin{center}
{\renewcommand{\arraystretch}{2}
\begin{tabular}{c|c|c}

 $M_c(\sf triv)$ & $M_c(\sf sgn)$ & $M_c(\tau_l)$ \\
\hline
 $1-dc$ & $1+dc$ & $1$ \\

\end{tabular}}
\end{center}
\item

\item Let $d=2m$. The lowest weights of standard
$H_{k_1,k_2}$-modules are given in the table below:

\begin{center}
{\renewcommand{\arraystretch}{2}
\begin{tabular}{c|c|c|c|c}
$M_{k_1,k_2}(\sf triv)$ & $M_{k_1,k_2}(\eps_1)$  &
$M_{k_1,k_2}(\sf sgn)$  & $M_{k_1,k_2}(\eps_2)$  &
$M_{k_1,k_2}(\tau_l)$ \\
\hline
 $1-R$ & $1-R'$ &
 $1+R$ & $1+R'$  & $1$\\
\end{tabular}}
\end{center}
where $R=m(k_1+k_2)$ and $R'=m(k_1-k_2)$.

\een
\end{proposition}
\begin{proof}
The  proposition follows from the formula for $h_0 (V)$ (see
section \ref{O}).
\end{proof}

\subsection{Some irreducible finite dimensional representations.}
\begin{proposition}{\rm (\cite{CE}, Theorem 2.3)} \label{fingor}
Let $W$ be a real reflection group, and $V \subset M_{\sf
k}=M_{\sf k}(\sf triv)$ be a $W$-subrepresentation of dimension
$d=\dim(\h)$ sitting in degree $r$ and consisting of singular
vectors. Let $I_{\sf k}$ be the ideal generated by $V$. Assume
that the quotient representation $Q_{\sf k}=M_{\sf k}/I_{\sf k}$
is finite dimensional. Then the representation $Q_{\sf k}$ is
irreducible. Its character is given by the formula
$$
\chi_{Q_{\sf k}}(g,t)=t^{h_0({\sf triv})} \frac{\det\bigl
|_V(1-gt^r)}{\det|_{\h^*}(1-gt)}.
$$
In particular, the dimension of $Q_{\sf k}$ is $r^d$.
\end{proposition}

\section{Proofs}
\begin{remark}
In all the classification theorem the character formulas follow
easily from the structure of standard modules.
\end{remark}

\subsection{Case of $d=2m+1$.}

\begin{proof}[Proof of theorem \ref{odd2dim}]
Let us describe the $W$-action on the set $U'$ of singular vectors
of positive degree.

{\it Claim $1$. There are no two-dimensional irreducible
subrepresentations in $U'$.}

By proposition \ref{lowesweights}, $\mathbf h$ acts by $1$ on all
irreducible two-dimensional representations of $W$ consisting of
singular vectors. However, the space on which $\mathbf h$ acts by
$1$, is exactly the space of vectors of degree $0$, which does not
intersect $U'$.

\medskip

{\it Claim $2$. For $c \ne n \pm l/d$ there are no one-dimensional
representations in $U'$.}

By symmetry reasons (see  remark \ref{symmetry}), it is enough to
prove that $U'$ does not contain the sign representation.

Suppose $U'$ contains a $W$-representation isomorphic to ${\sf
sgn}$. This representation generates an $H_c$-submodule isomorphic
to $M_c({\sf sgn})$ in $M_c(\tau_l)$. Therefore

$$
 \Hom_{H_c}(M_c({\sf sgn}), M_c(\tau_l)) \ne 0.
$$

Both $M_c({\sf sgn})$ and $M_c(\tau_l)$ are free over $\CC [\h]$.
By theorem \ref{equiv2},
 \beq \label{homineq}
\dim \Hom_{H_{\sf k}}(M_{c}({\sf sgn}),M_{c}(\tau_l)) \leqslant
\dim \Hom_{\mathcal H} ({\sf KZ}(M_c({\sf sgn})), {\sf
KZ}(M_c(\tau_l))),
 \eeq
and hence,

\beq \label{1dimhom}
 \Hom_{\mathcal H}({\sf KZ}(M_c({\sf sgn})), {\sf
KZ}(M_c(\tau_l)))\ne 0. \eeq

Let us compute  the dimension of $\Hom_{\mathcal
H\text{-mod}}({\sf KZ}(M_c({\sf sgn})), {\sf KZ}(M_c(\tau_l)))$.

Recall that in case of $W=I_2(d)$ the Hecke algebra $\mathcal
H(q)$ is generated by two elements $T_1$ and $T_2$ with relations
$$
(T_j-1)(T_j+q)=0, \quad \underbrace{T_1T_2 \dots
T_1}_{\mbox{\scriptsize $d$ factors}}=\underbrace{T_2T_1 \dots
T_2}_{\mbox{\scriptsize $d$ factors}}.
$$
In order for theorems \ref{equiv1} and \ref{equiv2} to hold, we
need $q=e^{-2\pi i c}$.

Representation ${\sf KZ}(M_c(\tau_l))$ of $\mathcal H(q)$
corresponds to $\tau_l$ via Tits' Deformation Theorem. If $0 < c
\ne n \pm l/d$, this representation is given explicitly by:
$$
T_1 \mapsto \left(
\begin{array}{cc}
1 & 0\\
c_l & -q
\end{array}\right)
; \quad T_2 \mapsto \left(
\begin{array}{cc}
-q & c_l'\\
0 & 1
\end{array}\right),
$$
where $c_l c_l'=q(2+\omega^l+ \omega^{-l})$ (see \cite{GP}, p.
268; note that in our case the sign of the generators is
different). Representation ${\sf KZ}(M_c({\sf sgn}))$ is given
explicitly by $T_1 \mapsto -q$ and $T_2 \mapsto -q$.

Hence, $
 \dim \bigl [\Hom_{\mathcal H}({\sf KZ}(M_c({\sf
 sgn})), {\sf KZ}(M_c(\tau_l)))\bigr ]= 0 $ for $0 < c
\ne n \pm l/d$.
 This contradicts
to (\ref{1dimhom}), and hence $U'$ does not contain
subrepresentations isomorphic to ${\sf sgn}$. Claim $2$ is proved.

\smallskip

From Claims $1$ and $2$ it follows that $U'=0$ for $r \not\equiv
\pm l \text{ $($mod $d)$}$. Thus, $U$ is the space of all vectors
of degree $0$ in $M_c(\tau_l)$. It is isomorphic to $\tau_l$ as a
representation of $W$, and spans $M_c(\tau_l)$. By proposition
\ref{irred}, $M_c(\tau_l)$ is irreducible for these values of $c$.
Hence, the last row of the table is proved.

\medskip

{\it Claim $3$. If $0<c=n \pm l/d$ for some $n \in \mathbb Z$,
then $U' \simeq {\sf sgn}$ as a $W$-representation. The
$H_c(W)$-submodule $I$, spanned by $U'$,  and  the quotient
$M_c(\tau_l)/I$ are both irreducible.}

By claim $1$, $U'$ does not contain two-dimensional
representations of $W$. By proposition \ref{lowesweights}, for
$c>0$ it does not contain the trivial representation either.
Hence, $U'$ is either $0$ or several copies of $\sf sgn$.

On the one hand, the multiplicity of $\sf sgn$  in $U'$ is at most
1 (specialize to one parameter results of \cite{Du1}, p.~82 and
use inequality (\ref{homineq})). On the other hand, according to
the theorem \ref{oddsingvect}, there is a singular copy of
$\tau_l$ in $M_c({\sf triv})$. An $H_c$-submodule $J$, generated
by this copy, is not isomorphic to $M_c(\tau_l)$, since $M_c({\sf
triv})$ is free of rank 1 over $\CC [\h]$ meanwhile $M_c(\tau_l)$
is free of rank 2. Hence $M_c(\tau_l)$ is reducible and $U' \ne
0$. This proves that $U' \simeq {\sf sgn}$.

The submodule $I$ generated by $U'$ is isomorphic to $M_c({\sf
sgn})$. By theorem \ref{irred} it is irreducible.

From proposition \ref{lowesweights} it follows that all the
singular vectors in $M_c(\tau_l)/I$ are images of singular vectors
in $M_c(\tau_l)$. From this fact and proposition \ref{irred} it
follows, that $M_c(\tau_l)/I$ is irreducible. This completes the
proof of claim $3$ and of the first row in the table.

The second row in the table follows via symmetry arguments.
\end{proof}

\medskip

\begin{proof}[ Proof of theorem \ref{odd1dim}]
\hspace{8cm}

{\it Claim  1. If $0< r \equiv \pm l$ (mod $d$), then $M_c({\sf
triv})$ contains a submodule isomorphic to $L_c(\tau_l)$, the
quotient by which is irreducible.}

By theorem \ref{oddsingvect}, $U' \subset M_c({\sf triv})$ is
isomorphic to $\tau_l$ as a representation of $W$, and consists of
vectors of degree $r$. Denote an $H_c$-submodule of $M_c({\sf
triv})$ generated by $U'$ by $I_c$. By definition, $I_c$ is a
module with lowest weight $\tau_l$. By theorem \ref{odd2dim},
there are only two such modules: $M_c(\tau_l)$ and $L_c(\tau_l)
\simeq M_c(\tau_l)/M_c({\sf sgn})$. Since there are no singular
vectors of type $\sf sgn$ in $I_c$ (see theorem
\ref{oddsingvect}), $I_c \simeq L_c(\tau_l)$.

For eigenvalues of $\mathbf h$ big enough, the dimensions of the
corresponding eigenspaces in $M_{c}({\sf triv})$ and $L_c(\tau_l)$
coincide. Hence, the quotient $M_c({\sf triv})/L_{c}(\tau_l)$ is
finite dimensional. The claim follows from proposition
\ref{fingor}.

\medskip

{\it Claim 2. If $0<c=n+\half$, then $M_c({\sf triv})$ contains a
submodule isomorphic to $L_c({\sf sgn })$, the quotient by which
is irreducible.}

 By theorem \ref{oddsingvect}, $U' \subset M_c({\sf triv})$ is
isomorphic to $\sf sgn$ as a $W$-representation, and is spanned by
a vector of degree $d(2n+1)$. A submodule generated by $U'$ is
isomorphic to $M_c({\sf sgn})$. By symmetry reasons and claim 1,
it is irreducible, i.e. isomorphic to $L_c({\sf sgn})$.

By the lowest weights reasons (proposition \ref{lowesweights}),
there are no singular vectors of positive degree in $Q_c=M_c({\sf
triv})/M_c({\sf sgn})$, and hence $Q_c$ is irreducible. Comparing
the dimensions of $\mathbf h$-eigenspaces one  can show, that
$Q_c$ is infinite dimensional.

\medskip

{\it Claim 3. $M_c({\sf triv})$ is irreducible for all other
values of $c$.}

By theorem \ref{oddsingvect}, there are no singular vectors in
$M_c({\sf triv})$ for all other values of $c$. Hence, this module
is irreducible.
\end{proof}

\subsection{Case of $d=2m$.}

\begin{proof}[Proof of theorem \ref{even2dim}] The proof is analogous to
that of theorem \ref{odd2dim}. We will describe the action of $W$
on the set $U'$ of singular vectors of positive degree.

{\it Claim 1. There are no two-dimensional irreducible
subrepresentations in $U'$.}

See proof of claim 1 of theorem \ref{odd2dim}.

\medskip

{\it Claim 2. If $(k_1,k_2) \notin E_r^+$ for all $r \equiv \pm l
 \text{ $($mod } 2m)$ and $(k_1,k_2) \notin E^-_{r'}$ for all $r'
\equiv m \pm l \text{ $($mod $2m)$}$, then $U'$ does not contain
any one-dimensional subrepresentations of $W$.}

By symmetry reasons (remark \ref{symmetry}), it is enough to prove
the following fact.

{\it Claim 2 $'$.
 If $(k_1,k_2) \notin E_r^+$ for all $r \equiv
\pm l \text{
 $($mod $2m)$}$, then $U'$ does not contain any $W$-subrepresentations
isomorphic to $\sf sgn$.}

Suppose $U'$ contains a $W$-representation isomorphic to ${\sf
sgn}$. This representation generates an $H_{k_1,k_2}$-submodule of
$M_{k_1,k_2}(\tau_l)$, which is isomorphic to $M_{k_1,k_2}({\sf
sgn})$. Therefore
 $$
 \Hom_{H_{k_1,k_2}}(M_{k_1,k_2}({\sf sgn}), M_{k_1,k_2}(\tau_l)) \ne 0.
 $$

Both $M_{k_1,k_2}({\sf sgn})$ and $M_{k_1,k_2}(\tau_l)$ are free
over $\CC [\h]$. By theorem \ref{equiv2},

\beq \label{homineq2} \text{ $
\begin{array}{l}
\dim \Hom_{H_{\sf k}}(M_{k_1,k_2}({\sf sgn}),M_{k_1,k_2}(\tau_l))
\leqslant
\vspace{5pt}\\
\hspace{4cm} \leqslant \dim \Hom_{\mathcal H} ({\sf
KZ}(M_{k_1,k_2}({\sf sgn})), {\sf KZ}(M_{k_1,k_2}(\tau_l))),\\
\end{array}$}
\eeq and hence,
 \beq \label{2dimhom}
 \Hom_{\mathcal H\text{-mod}}({\sf
KZ}(M_{k_1,k_2}({\sf sgn})), {\sf KZ}(M_{k_1,k_2}(\tau_l)))\ne 0.
 \eeq

Let us compute  the dimension of $\Hom_{\mathcal
H\text{-mod}}({\sf KZ}(M_{k_1,k_2}({\sf sgn})), {\sf
KZ}(M_{k_1,k_2}(\tau_l)))$.

Recall that in case of $W=I_2(2m)$ the Hecke algebra $\mathcal
H(q_1,q_2)$ is spanned by two elements $T_1$ and $T_2$ with
relations
$$
(T_j-1)(T_j+q_j)=0, \quad \underbrace{T_1T_2 \dots
T_2}_{\mbox{{\scriptsize $2m$ multiples}}}=\underbrace{T_2T_1
\dots T_1}_{\mbox{{\scriptsize $2m$ multiples}} },
$$

In order for theorems \ref{equiv1} and \ref{equiv2} to hold, we
need $q_j=e^{-2\pi i k_j}$.

The representation ${\sf KZ}(M_{k_1,k_2}(\tau_l))$ of $\mathcal
H(q_1,q_2)$ corresponds to $\tau_l$ via Tits' Deformation Theorem.
For $(k_1,k_2) \notin E_r^+$ for all $0<r \equiv \pm l \text{ (mod
$2m$)}$, this representation is given explicitly by:
$$
T_1 \mapsto \left(
\begin{array}{cc}
1 & 0\\
c_l & -q_1
\end{array}\right)
; \quad T_2 \mapsto \left(
\begin{array}{cc}
-q_2 & c_l'\\
0 & 1
\end{array}\right),
$$
where $c_l c_l'=q_1+q_2 + \sqrt{q_1q_2}(\omega^l+ \omega^{-l})$
and $\sqrt{q_1q_2}$ is fixed to be $e^{-\pi i (k_1+k_2)}$ (see
\cite{GP}, p.~268; note that in our the signs of the generators
are different). Using the notations of \cite{Du1}, these are
formulas for $-e^{2\pi i \alpha}M(\gamma_0)$ and $-e^{2\pi i
\beta}M(\gamma_1)$ in the basis $\{f^{(3)}, f^{(1)}\}$.

The representation ${\sf KZ}(M_c({\sf sgn}))$ is given explicitly
by $T_1 \mapsto -q_1$ and $T_2 \mapsto -q_2$.

Hence, $ \dim \Hom_{\mathcal H}({\sf KZ}(M_{k_1,k_2}({\sf sgn})),
{\sf KZ}(M_{k_1,k_2}(\tau_l))=0$, if $(k_1,k_2) \notin E_r^+$ for
all $0<r \equiv \pm l \text{ (mod $2m$)}$. This contradicts to
(\ref{2dimhom}), and hence $U'$ does not contain
subrepresentations isomorphic to ${\sf sgn}$. Claim $2$ is proved.

From Claims $1$ and $2$ it follows that $U'=0$, when $(k_1,k_2)
\notin E^+_r$ and $(k_1,k_2) \notin E^-_{r'}$ for all $r \equiv
\pm l \text{ $($mod $2m)$}$ and all $r' \equiv m\pm l \text{
$($mod $2m$)}$ . By proposition \ref{irred}, $M_{k_1,k_2}(\tau_l)$
is irreducible for these values of $(k_1,k_2)$. This proves part
{\bf (iv)} of the theorem.

\medskip

{\it Claim 3. If $(k_1,k_2) \in E^+_r$ for some $0<r \equiv \pm l
\text{ $($mod $2m)$}$ and $(k_1,k_2) \notin E^-_{r'}$ for all $r'
\equiv m \pm l \text{ $($mod $2m)$ }$, then $U' \simeq {\sf sgn}$
as a representation of $W$. The $H_{k_1,k_2}(W)$-submodule
$I_{k_1,k_2}$, spanned by $U'$, and the quotient
$M_{k_1,k_2}(\tau_l)/I_{k_1,k_2}$ are both irreducible.}

Let us prove that $U'$ is a direct sum of representations,
isomorphic to $\sf sgn$. By claim $1$, $U'$ does not contain any
two-dimensional representations of $W$. By proposition
\ref{lowesweights}, for $r>0$ it does not contain the trivial
representation either. By symmetry reasons (remark \ref{symmetry})
and claim $2'$, the $W$-representation $U'$ does not contain any
subrepresentations of type $\eps_1$ and $\eps_2$. Hence $U'$ is
either $0$ or several copies of $\sf sgn$.

Let us prove that the multiplicity of $\sf sgn$ is exactly $1$. On
the one hand, the multiplicity of $\sf sgn$ in $U'$ is at most 1
(see inequality (\ref{homineq2}) and \cite{Du1}, p. 82).
 On the
other hand, according to
 theorem \ref{oddsingvect}, there is a singular copy
of $\tau_l$ in $M_{k_1,k_2}({\sf triv})$. An
$H_{k_1,k_2}$-submodule $J$, generated by this copy, is not
isomorphic to $M_{k_1,k_2}(\tau_l)$, since $M_{k_1,k_2}({\sf
triv})$ is free of rank 1 over $\CC [\h]$ meanwhile
$M_{k_1,k_2}(\tau_l)$ is free of rank 2. Hence,
$M_{k_1,k_2}(\tau_l)$ is reducible and $U' \ne 0$. This proves
that $U' \simeq {\sf sgn}$.

The submodule $I_{k_1,k_2}$, generated by $U'$, is isomorphic to
$M_{k_1,k_2}({\sf sgn})$. By theorem \ref{irred}, it is
irreducible.

Denote by $Q_{k_1,k_2}$ the quotient
$M_{k_1,k_2}(\tau_l)/I_{k_1,k_2}$. We have already proved that all
the singular vectors of positive degree in $Q_{k_1,k_2}$ are
secondary, i.e. singular only "modulo $I_{k_1,k_2}$". By
proposition \ref{lowesweights}, there are two possibilities for
secondary singular vectors:

\begin{itemize}
\item{if $(k_1,k_2) \in {E_{r'}^-}$ for some $r'>0$ and $k_2$ is negative,
then $Q_{k_1,k_2}$ might contain a singular vector of type
$\eps_2$;}
\item{if $(k_1,k_2) \in {E_{r'}^-}$ for some $r'<0$ and $k_1$ is negative,
then $Q_{k_1,k_2}$ might contain a singular vector of type
$\eps_1$.}
\end{itemize}

Suppose $Q_{k_1,k_2}$ contains a singular vector of type $\eps_2$.
The $H_{k_1,k_2}$-submodule $P_{k_1,k_2}$, generated by this copy
is isomorphic to $L_{k_1,k_2}(\eps_2) \simeq M_{k_1,k_2}(\eps_2)$.
Denote by $N$ the preimage of $P_{k_1,k_2}$ in
$M_{k_1,k_2}(\tau_l)$. Since $P_{k_1,k_2}$ is generated by
secondary vector, the following exact sequence does not split.
$$
0 \rightarrow M_{k_1,k_2}({\sf sgn}) \rightarrow N \rightarrow
M_{k_1,k_2}(\eps_2) \rightarrow 0
$$
Apply ${\sf KZ}$-functor to this sequence. The representation
${\sf KZ}(M_{k_1,k_2}({\sf sgn}))$ is given by $T_j \mapsto -q_j$.
The representation ${\sf KZ}(M_{k_1,k_2})(\eps_1)$ is given by
$T_1 \mapsto -q_1$ and $T_2 \mapsto 1$. To find ${\sf KZ}(N)$
consider the following exact sequence
$$
0 \rightarrow N \rightarrow M_{k_1,k_2}(\tau_l) \rightarrow J
\rightarrow 0.
$$
The cokernel $J$ is supported on $\h \backslash \h_{reg}$, and
hence ${\sf KZ}(J)=0$. Applying $\sf KZ$-functor to the sequence
above, we will get that ${\sf KZ}(N) \simeq {\sf
KZ}(M_{k_1,k_2}(\tau_l))$, since the functor is exact. Hence,
${\sf KZ}(N)$ is given explicitly by
$$
T_1 \mapsto \left(
\begin{array}{cc}
1 & 0\\
c_l & -q_1
\end{array}\right)
; \quad T_2 \mapsto \left(
\begin{array}{cc}
-q_2 & c_l'\\
0 & 1
\end{array}\right).
$$
Since $\sf KZ$-functor is exact, ${\sf KZ}(N)$ is a nontrivial
extension of ${\sf KZ}(M_{k_1,k_2}(\eps_2))$ by ${\sf
KZ}(M_{k_1,k_2}({\sf sgn}))$. It is impossible, since in every
such extension $\Tr (T_1) =-2 q_1$ meanwhile the trace of $T_1$ in
${\sf KZ}(N)$ is $1-q_1$, which is not equal to $-2q_1$ for $r'
\not\equiv m \pm l \text{ (mod $2m$)}$. Hence, $Q_{k_1,k_2}$ does
not contain any submodules and is irreducible.

The case of singular vector of type $\eps_1$ is proved
analogously. Claim 3 is proved.


\medskip

{\it Claim 4. Let $k_1 \geqslant k_2 >0$ be such that $(k_1,k_2)
\in E^+_r \cap E^-_{r'}$ for some $r \equiv \pm l \text{ $($mod
$2m)$}$ and some $r' \equiv m \pm l \text{ $($mod $2m)$}$. Then
$M_{k_1,k_2}(\tau_l)$ contains a submodule $I_{k_1,k_2}$,
isomorphic to $M_{k_1,k_2}({\sf sgn})$, and a submodule
$\widetilde{I}_{k_1,k_2}$, isomorphic to $M_{k_1,k_2}(\eps_2)$. If
$k_2 \in \mathbb Z +1/2$, then $\widetilde{I}_{k_1,k_2} \supset
I_{k_1,k_2}$, otherwise $\widetilde{I}_{k_1,k_2} \cap
I_{k_1,k_2}=0$.}

\smallskip

By claim $3$, $M_{k_1,k_2}(\tau_l)$ contains a submodule
$I_{k_1,k_2}$, isomorphic to $M_{k_1,k_2}({\sf sgn})$, when
$(k_1,k_2)$ is a generic point of $E^+_r$, and it contains a
submodule $\widetilde{I}_{k_1,k_2}$, isomorphic to
$M_{k_1,k_2}(\eps_2)$, when $(k_1,k_2)$ is a generic point of
$E^-_{r'}$. By continuity for $(k_1,k_2) \in E^+_r \cap E^-_{r'}$,
the standard module $M_{k_1,k_2}(\tau_l)$ contains both of them.

Suppose $I_{k_1, k_2}  \cap \widetilde{I}_{k_1, k_2} \ne 0$. Let
$f \in I_{k_1, k_2}  \cap \widetilde{I}_{k_1, k_2} $ be a vector
of the minimal possible degree. Then $f$ is a singular vector in
both $I_{k_1, k_2}$ and $\widetilde{I}_{k_1, k_2}$. The only
singular vector in $I_{k_1,k_2}$ is the one which generates this
submodule. Indeed, by symmetry arguments, the structure of
$M_{k_1, k_2}(\sf sgn)$ is the same as the structure of $M_{-k_1,
-k_2}(\sf triv)$. By proposition \ref{evensingvect}, the unique
singular vector contained in $M_{-k_1, -k_2}(\sf triv)$ generates
the whole module. Hence, if intersection of $I_{k_1, k_2}$ and
$\widetilde{I}_{k_1, k_2}$ is nonzero, then $\widetilde{I}_{k_1,
k_2} \supset I_{k_1, k_2}$.

Suppose this intersection is nonzero. Then in
$M_{k_1,k_2}(\eps_2)$ there is a singular vector, on which
$I_2(2m)$ acts by $\sf sgn$. By symmetry reasons, this happens
only for those $(k_1,k_2)$, for which $M_{-k_1, k_2}(\sf triv)$
contains a submodule isomorphic to $M_{-k_1, k_2}(\eps_1)$. By
proposition \ref{evensingvect}, this is the case  if $k_2=i+\half$
for some positive integer $i$. This proves claim 4.

\medskip

{\it Claim 5. Let $(k_1,k_2)$ be as in claim 4. Then the quotient
$Q_{k_1,k_2}$ of $M_{k_1,k_2}(\tau_l)$ by the submodule generated
by $I_{k_1,k_2}$ and $\widetilde{I}_{k_1,k_2}$ is irreducible.}

\smallskip

By proposition \ref{lowesweights}, the only possible singular
vectors in $Q_{k_1,k_2}$ are of type ${\sf sgn}$ and degree $r$,
or of type $\eps_2$ and degree $r'$.

From inequality (\ref{homineq2}) and \cite{Du1}, p. 84, it follows
that the set of singular vectors in $M_{k_1,k_2}(\tau_l)$ contains
only one copy of $\sf sgn$ and only one copy of $\eps_2$, i.e. is
isomorphic to $\tau_l \oplus \eps_2 \oplus {\sf sgn}$ as a
$W$-representation. Hence, all the singular vectors in
$Q_{k_1,k_2}$, except the one, which generates this module, are
secondary, i.e. singular only "modulo $\widetilde{I}_{k_1,k_2}$
and $I_{k_1,k_2}$".

For the values of parameters that we consider, the lowest possible
degree of a vector in the ideal generated by
$\widetilde{I}_{k_1,k_2}$ and $I_{k_1,k_2}$ is $r'$. Hence, a
vector of type $\eps_2$ can not be secondary. Therefore, if
$Q_{k_1,k_2}$ is reducible, then it contains a singular vector of
type $\sf sgn$. By theorem \ref{lowesweights}, the
$H_{k_1,k_2}$-submodule $P_{k_1,k_2}$, generated by this copy, is
isomorphic to $L_{k_1,k_2} ({\sf sgn}) \simeq M_{k_1,k_2}({\sf
sgn})$.

For $k_2 \notin \mathbb Z + 1/2$ the intersection of
$\widetilde{I}_{k_1,k_2}$ and $I_{k_1,k_2}$ is zero. Suppose that
$Q_{k_1,k_2}$ has a submodule $P_{k_1,k_2}$. Then
$Q_{k_1,k_2}/P_{k_1,k_2}$, and hence $L_{k_1,k_2}(\tau_l)$, is
finite dimensional. But it is known (see \cite{De} remark 5.7),
that $L_{k_1,k_2}(\tau_l)$ is always infinite dimensional. Hence,
$Q_{k_1,k_2}$ is irreducible.

For $k_2 \in \mathbb Z + 1/2$ we have the following sequence of
modules
$$
M_{k_1,k_2}(\tau_l)\supset M_{k_1,k_2}(\eps_2) \supset
M_{k_1,k_2}({\sf sgn}),
$$
and $Q_{k_1,k_2} = M_{k_1,k_2}(\tau_l) / M_{k_1,k_2}(\eps_2)$.
Proof of the irreducibility of $Q_{k_1,k_2}$ is analogous to the
proof of the irreducibility of the quotient in claim 3. Claim 5 is
proved.

\medskip

Part {\bf (iii)} of the theorem for $k_1 \geqslant k_2 >0$ follows
from claims 4 and 5. For all the other values of parameters the
proof is analogous.
\end{proof}

\begin{proof}[Proof of theorem \ref{even1dim}]

{\bf (i)} By theorem \ref{evensingvect},  the space $U' \subset
M_{k_1, k_2}({\sf triv})$ is isomorphic to $\tau_l$ as a
representation of $W$ and sits in degree $r$. The submodule
$I_{k_1,k_2}$, spanned by $U'$, is isomorphic to
$L_{k_1,k_2}(\tau_l)$ since it does not contain any singular
vectors except the one, by which it is generated. By  parts {\bf
(i)} and {\bf (iii)} of theorem \ref{even2dim},
$L_{k_1,k_2}(\tau_l) \simeq M_{k_1,k_2}(\tau_l)/M_{k_1,k_2}({\sf
sgn})$.

For eigenvalues of $\mathbf h$ big enough, the dimensions of the
corresponding eigenspaces in $M_{k_1, k_2}(\sf triv)$ and $L_{k_1,
k_2}(\tau_l)$ coincide. Hence, the quotient $M_{k_1,k_2}({\sf
triv})/L_{k_1,k_2}(\tau_l)$ is finite dimensional. The statement
follows from proposition \ref{fingor}.

{\bf (ii)} By theorem \ref{evensingvect}, the space $U' \subset
M_{k_1, k_2}(\sf triv)$ is isomorphic to $\eps_2$ as a
$W$-representation. It is easy to see that it generates a
submodule isomorphic to $M_{k_1, k_2}(\eps_2)$.

Consider the quotient  $Q_{k_1,k_2} =M_{k_1,k_2}({\sf
triv})/M_{k_1,k_2}(\eps_2)$. From the proposition
\ref{lowesweights} it follows that $Q_{k_1,k_2}$ does not contain
any singular vectors and hence is irreducible. Comparing
dimensions of eigenspaces of $\mathbf h$ in $M_{k_1,k_2}(\eps_2)$
and in $M_{k_1,k_2}({\sf triv})$, one can show that $Q_{k_1,k_2}$
is infinite dimensional.

{\bf (iii)} Analogous to part ${\bf (ii)}$.

{\bf (iv)} Consider the case of $(k_1,k_2)\in E^+_r \cap L_i^1$
with $k_1<k_2$. The other case is proved analogously.

 By theorem
\ref{evensingvect}, the space $U' \subset M_{k_1,k_2}({\sf triv})$
is isomorphic to $\tau_l \oplus \eps_2$ as a $W$-representation.

A submodule $I_{k_1,k_2}$ of $M_{k_1,k_2}({\sf triv})$, generated
by a singular vector of type $\eps_2$, is isomorphic to
$M_{k_1,k_2}(\eps_2)$. By symmetry reasons, $M_{k_1,k_2}(\eps_2)$
has the same structure as $M_{-k_1,k_2}({\sf triv})$ and hence, by
part {\bf (i)}, it contains a singular copy of $\tau_{l}$, which
generates a submodule isomorphic to $L_{k_1,k_2}(\tau_l)$. So we
have the following picture:
$$
M_{k_1,k_2}({\sf triv}) \supset M_{k_1,k_2}(\eps_2) \supset
L_{k_1,k_2}(\tau_l).
$$

By symmetry arguments,  since $M_{-k_1,k_2}({\sf
triv})/L_{-k_1,k_2}(\tau_{m-l})$ is irreducible, so is the
quotient $M_{k_1,k_2}(\eps_2)/L_{k_1,k_2}(\tau_l)$.

Consider the quotient $Q_{k_1,k_2}=M_{k_1,k_2}({\sf
triv})/M_{k_1,k_2}(\eps_2)$. Comparing the dimensions of $\mathbf
h$-eigenspaces in $M_{k_1,k_2} ({\sf triv})$ and $M_{k_1,
k_2}(\eps_2)$, one gets that $Q_{k_1,k_2}$ is infinite
dimensional.

Since  $Q_{k_1,k_2}$ has the lowest weight $\sf triv$, its
irreducible quotient is isomorphic to $L_{k_1,k_2}({\sf triv})$.
For all $(k_1,k_2)$ the module $L_{k_1,k_2}({\sf triv})$ is
defined as a quotient of $M_{k_1,k_2}({\sf triv})$ by the kernel
of Shapovalov form. By part {\bf (i)}, for generic $(k_1,k_2) \in
E_r^+$ the irreducible quotient has dimension $r^2$ and hence
$\dim L_{k_1,k_2}({\sf triv}) \leqslant r^2$.

Hence $Q_{k_1,k_2}$ is reducible and contains singular vectors.
Note that these vectors are secondary, i.e. singular only ``modulo
$M_{k_1,k_2}(\eps_2)$". By the table of lowest $\mathbf h$-weights
(proposition \ref{lowesweights}), there exit the following
possibilities for singular vectors: vectors of type $\tau_l$
sitting in degree $r$, vectors of type $\eps_1$ sitting in degree
$2mk_2 $ or vectors of type $\sf sgn$ in degree $2r$.

The sequence of the dimensions of $\mathbf h$-eigenspaces in
$Q_{k_1,k_2}$ is
 \beq \label{hdims}
 (1,2,3, \dots, 2k_1m, 2k_1m, 2k_1m, \dots).
 \eeq
Since $L_{k_1,k_2}({\sf triv})$ is a finite dimensional
$\mathfrak{sl}_2$-module, the sequence of dimensions of its
$\mathbf h$-eigenspaces is finite and symmetric with respect to
the $0$-eigenspace. From the table of lowest $\mathbf h$-weights
(proposition \ref{lowesweights}) it follows that $\mathbf h$ acts
by $0$ on polynomials of degree $r-1$. In order to make the
sequence (\ref{hdims}) symmetric, $Q_{k_1,k_2}$ has to contain a
unique up to proportionality singular vector, which is of type
$\eps_1$ and sits in degree $2mk_2$. This vector generates a
submodule $P_{k_1,k_2}$ isomorphic to $L_{k_1,k_2}(\eps_1)$. From
the sequence of dimensions of $\mathbf h$-eigenspaces in
$Q_{k_1,k_2}/P_{k_1,k_2}$ it follows that this quotient is is
irreducible.


{\bf (v)} Consider the case $(k_1,k_2)\in E_r^+ \cap L_i^1$. By
theorem \ref{evensingvect}, there exists a singular vector of type
$\eps_2$ in $M_{k_1,k_2}({\sf triv})$. It generates a submodule
isomorphic to $M_{k_1,k_2}(\eps_2)$. This submodule is irreducible
since, by symmetry arguments, it has the same structure as
$M_{-k_1,k_2}({\sf triv})$.

Consider the quotient $Q_{k_1,k_2}=M_{k_1,k_2}({\sf
triv})/M_{k_1,k_2}(\eps_2)$. By theorem \ref{evensingvect}, there
is a singular copy of $\tau_l$ in $Q_{k_1,k_2}$. Denote by
$P_{k_1,k_2}$ the submodule generated by this copy. By definition,
this is a module of lowest weight $\tau_l$. By theorem
\ref{even2dim}, for these values of parameters there are three
such modules: $M_{k_1,k_2}(\tau_l)$,
$M_{k_1,k_2}(\tau_l)/M_{k_1,k_2}({\sf sgn})$ and $L_{k_1,k_2}({\sf
triv}) \simeq M_{k_1,k_2}({\sf triv})/M_{k_1,k_2}(\eps_2) $. Since
there are no singular vectors of type ${\sf sgn}$ and $\eps_2$ in
$Q_{k_1,k_2}$, the submodule $P_{k_1,k_2}$ is isomorphic to
$L_{k_1,k_2}(\tau_l)$.

Using the symmetry of the sequence of the dimensions of $\mathbf
h$-eigenspaces, one can show that $Q_{k_1,k_2}/P_{k_1,k_2}$ is
irreducible.

{\bf (vi)} By theorem \ref{evensingvect}, in $M_{k_1,k_2}(\sf
triv)$ there are singular vectors of types $\eps_1$ and $\eps_2$.
They generates submodules isomorphic to  $M_{k_1,k_2}(\eps_1)$ and
$ M_{k_1,k_2}(\eps_2)$ respectively. These two modules intersect
each other by a submodule isomorphic to $M_{k_1,k_2}({\sf sgn})$
and generated by a singular vector. By symmetry reasons, the
quotients $M_{k_1,k_2}(\eps_1)/M_{k_1,k_2}({\sf sgn})$ and
$M_{k_1,k_2}(\eps_2)/M_{k_1,k_2}({\sf sgn})$ are irreducible. The
quotient of $M_{k_1,k_2}({\sf triv})$ by the union of
$M_{k_1,k_2}(\eps_1)$ and $M_{k_1,k_2}(\eps_2)$ is irreducible for
the same reasons as in parts {\bf (iv)} and {\bf (v)}.

{\bf (viii)} By theorem \ref{evensingvect}, the space $U'$ is
empty and hence, by proposition \ref{irred} $M_{k_1,k_2}({\sf
triv})$ is irreducible.
\end{proof}

{\bf Acknowledgements.} I am very grateful to P.~Etingof for
suggesting me this problem and for many useful conversations. I
 would also like to thank I.~Scherbak and S.~Chmutov for many
comments.

\end{document}